\documentclass[A4,14pt]{amsart}

\usepackage{inputenc} %[latin1]
\usepackage[T1]{fontenc}
\usepackage{amssymb,amsfonts,amsmath,amsthm}
\newtheorem{Theorem}{Theorem}
\newtheorem{Corollary}[Theorem]{Corollary}

\newtheorem{theorem}{Theorem}[section]
\newtheorem{proposition}[theorem]{Proposition}
\newtheorem{corollary}[theorem]{Corollary}
\newtheorem{lemma}[theorem]{Lemma}

\newtheorem{definition}[theorem]{Definition}
\newtheorem{problem}[theorem]{Problem}

\newtheorem{fact}[theorem]{Fact}
\newtheorem{claim}[theorem]{Claim}

\newtheorem{property}[theorem]{Property}
\newtheorem{example}[theorem]{Example}

\newcommand{\MM}{{\sf MM}}
\newcommand{\PFA}{{\sf PFA}}

\newcommand{\ZFC}{{\sf ZFC}}
\newcommand{\SCH}{\sf SCH}

\newcommand{\CP}{\sf CP}

\newcommand{\w}{\omega}
\renewcommand{\a}{\alpha}

\newcommand{\cf}{{\mbox{\tt cof}}}
\newcommand{\otp}{{\mbox{\tt otp}}}
\newcommand{\D}{{\mathcal D}}
\newcommand{\E}{{\mathcal E}}
\newcommand{\I}{{\mathcal I}}

\newcommand{\C}{\mathcal C}
\newcommand{\A}{\mathcal A}

\newcommand{\pcf}{{\sf pcf}}
\newcommand{\la}{{\langle}}
\newcommand{\ra}{{\rangle}}

\author{Assaf Sharon}
\author{Matteo Viale}
\thanks{The second author wishes to thank Boban Veli\v{c}kovi\'c for several useful hints and comments on previous drafts. In particular the results in subsection \ref{subsecWDR} are due to him.}
\title{Some consequences of reflection on the approachability ideal}

\begin{document}

\begin{abstract}
We study the approachability ideal $\I[\kappa^+]$ in the context of 
large cardinals properties of the regular cardinals below a singular $\kappa$. As a guiding example consider the
approachability ideal $\I[\aleph_{\w+1}]$ assuming that $\aleph_\w$ is strong limit. 
In this case we obtain that club many points in $\aleph_{\w+1}$ of cofinality $\aleph_n$ for some $n>1$ 
are approachable assuming the joint reflection of countable families of
stationary subsets of $\aleph_n$. This reflection principle
holds under $\MM$ for all $n>1$ and for each $n>1$ is equiconsistent with $\aleph_n$ being weakly compact in
$L$. This characterizes the structure of the approachability ideal $\I[\aleph_{\w+1}]$ in models of $\MM$.
We also apply our result to show that the Chang conjecture $(\aleph_{\w+1},\aleph_\w)\twoheadrightarrow(\aleph_2,\aleph_1)$ fails in models of $\MM$.
\end{abstract}

\maketitle
%%%%%%%%%%%%%%%%%%%%%%%%%%%%%%%%%%%%%%%%%%%%%%%%%%%%%%%%%%%%%%%%%%%%%%%%%%%%%%%%%%%%%%%%%%%%%%%%%%%%%%%%%%%%%%%%%%%%%%%%%%%%%%%%%%%%%%%%%%%%

\section{The approachability ideal}
In the course of development of the $\pcf$-theory of possible cofinalities Shelah
has introduced several interesting stationary sets on the successor of a singular
cardinal\footnote{\cite{eisHST} is our main reference source.}.
Among these are the sets of approachable and weakly approachable points in
$\kappa^+$, where $\kappa$ is a singular cardinal. 
Given $\A=\{a_\a:\a<\kappa^+\}\subseteq[\kappa^+]^{<\kappa}$, $\delta$ is weakly
approachable with respect to $\A$ if there is $H$ unbounded in $\delta$ of minimal
order type such that $\{H\cap\gamma:\gamma<\delta\}$ is covered\footnote{I.e.: for every $\gamma<\delta$ there is $\a<\delta$ such that $H\cap\gamma\subseteq a_\a$.} by $\{a_\a:\a<\delta\}$ and
$\delta$ is approachable with respect to $\A$ if there is $H$ unbounded in $\delta$
of minimal order type such that $\{H\cap\gamma:\gamma<\delta\}\subseteq\{a_\a:\a<\delta\}$.

\begin{definition} Let $\kappa$ be a singular cardinal. $S$ is (weakly) approachable
if there is a sequence $\A=\{a_\a:\a<\kappa^+\}\subseteq[\kappa^+]^{<\kappa}$ and a
club $C$ such that $\delta$ is (weakly) approachable with respect to $\a$
for all $\delta\in S\cap C$. $\I[\kappa^+]$ is the ideal generated by approachable
sets, $\I[\kappa^+,\kappa]$ is the ideal generated by weakly approachable sets.
\end{definition}

\noindent It is clear that $\I[\kappa^+]\subseteq\I[\kappa^+,\kappa]$.
For many of the known applications of approachability, it is irrelevant
whether we concentrate on the notion of weak approachability or on the apparently
stronger notion of approachability. Moreover in the case that $\kappa$ is strong
limit and singular $\I[\kappa^+]=\I[\kappa^+,\kappa]$ (section 3.4 and proposition 3.23 of
\cite{eisHST}). For this reason we feel free to concentrate our attention on the
notion of weak approachability which applies to a more general context.
It is rather easy to show that $\I[\kappa^+,\kappa]$ is a normal $\kappa^+$-closed
ideal which extends the non-stationary ideal.
A main result of Shelah is that there is a stationary set in $\I[\kappa^+]$
for any singular cardinal $\kappa$ (theorem 3.18 \cite{eisHST}).
There are several applications of this ideal to
the combinatorics of singular cardinals, we remind the reader one of them and refer
him to section 3 of \cite{eisHST} for a detailed account: the extent of this ideal
can be used to size the large cardinal properties of $\kappa$. $\I[\kappa^+,\kappa]$
is trivial unless the cardinals below $\kappa^+$ have very strong combinatorial
properties (in the range of supercompactness). Thus for example if
square at $\kappa$ holds $\I[\kappa^+]=\I[\kappa^+,\kappa]=P(\kappa^+)$ (theorem 3.13 of
\cite{eisHST}). On the other hand if $\lambda$ is strongly compact and $\kappa>\lambda$ is
singular of cofinality $\theta<\lambda$ then there is a stationary subset of
$\kappa^+$ of points of cofinality less than $\lambda$ which is not in
$\I[\kappa^+,\kappa]$ (Shelah, theorem 3.20 of \cite{eisHST}). In the same spirit if
$\MM$ holds there is a stationary set of points of cofinality $\aleph_1$ which is
not in $\I[\aleph_{\w+1},\aleph_\w]$ (Magidor, unpublished). It is also consistent\footnote{See for example
\cite{gitshaFAP} where this is achieved in the presence of a very good scale on
$\prod_{\a<\w^2}\aleph_\a$.}
that for unboundedly many $\a<\w^2$ there is a stationary set of points of
cofinality $\aleph_\a$ not in $\I[\aleph_{\w^2+1}]$. It is an open problem whether it is consistent that
there is a stationary set on $\aleph_{\w+1}$ concentrating on cofinalities larger
than $\aleph_1$ and not in $\I[\aleph_{\w+1}]$ (see for example the introduction of
\cite{formagVWS} or the end of section 3.5 in \cite{eisHST}). We will give a partial
answer to this question showing that this is not the case in models of $\MM$. Our results
have broader consequences and give serious constraints to the
possible scenarios where this problem may have a positive solution.
We briefly introduce some relevant concepts in our analysis. 
$S^\lambda_{\theta}$ denote the subset of $\lambda$ of points of cofinality
$\theta$. A stationary subset of $\lambda$ reflects on $\a$ if it intersects all the
closed and unbounded subsets of $\a$.

\begin{definition}
Let $\theta<\lambda$ be regular cardinals.

\smallskip

$R(\lambda,\theta)$ holds for infinite regular cardinals $\theta<\lambda$ if there is $S$
stationary subset of $\lambda$ such that for all families
$\{S_i:i<\theta\}$ of stationary subsets of $S$ there is $\delta<\lambda$ such that
$S_i$ reflects on $\delta$ for all $i$.

\smallskip

$R^*(\lambda)$ holds if if there is $S$ stationary subset of $\lambda$ such that for
all families
$\{S_i:i<\lambda\}$ of stationary subsets of $S$ there is $\delta<\lambda$ such that
$S_i$ reflects on $\delta$ for all $i<\delta$.
\end{definition}

\noindent It is clear that $R^*(\lambda)$ implies $R(\lambda,\zeta)$ which implies $R(\lambda,\theta)$ for all $\theta\leq\zeta<\lambda$.  
Moreover it is not hard to realize $R^*(\lambda)$ and $R(\lambda,\theta)$ and we will substantiate this in section \ref{secRP}. We now state one of our main result which gives rightaway a clear picture of what we are aiming to.
Given regular cardinals $\theta<\lambda$, $\lambda$ is $\theta$-inaccessible if $\zeta^\theta<\lambda$ for all $\zeta<\lambda$.

\begin{Theorem}\label{mainth}
Assume: 
\begin{itemize}
\item $\kappa$ is singular of cofinality $\theta$ and $\nu=\kappa^+$,
\item $\lambda<\kappa$ is either $\theta$-inaccessible or in $[\theta^+,\theta^{+\w})$,
\item $R(\lambda,\theta)$ holds.
\end{itemize}
Then $S_\lambda^{\nu}\in\I[\nu,\kappa]$.
\end{Theorem}

Immediate applications of theorem \ref{mainth} are the following:

\begin{Corollary} Assume $\lambda$ is weakly compact $\kappa>\lambda$ is singular cofinality $\theta<\lambda$ and $\nu=\kappa^+$.
Then $S^{\nu}_\lambda\in\I[\nu,\kappa]$.
\end{Corollary}
\begin{proof} $\lambda$ is $\theta$-inaccessible and satisfy $R(\lambda,\theta)$ (see fact \ref{fareflWC}). Now apply theorem \ref{mainth}.
\end{proof}

The reflection hypothesis of the main theorem holds in models of strong forcing axioms, for example we can prove:
\begin{Corollary} Assume Martin's maximum $\MM$ holds. Then club
many points in $S^{\aleph_{\w+1}}_{>{\aleph_1}}$ are approachable.
\end{Corollary}
\begin{proof} $\MM$ implies $R(\aleph_n,\aleph_1)$ holds as witnessed by $S^{\aleph_n}_\w$ for all $n>1$ (see \cite{formagsheMM}). Now apply theorem \ref{mainth}.
\end{proof}

We will also be able to obtain by a slight variation of the proof of theorem \ref{mainth}:
\begin{Theorem} 
Assume the proper forcing axiom $\PFA$. Then club
many points in $S^{\aleph_{\w+1}}_{>{\aleph_2}}$ are approachable.
\end{Theorem}

\begin{proof} By theorem \ref{mainth2} and theorem \ref{thPFACP}.
\end{proof}

Finally in section \ref{secCC} we will apply these results to the study of Chang conjecture and prove for example:
\begin{Theorem} 
Assume $R(\aleph_2,\aleph_0)$. Then $(\aleph_{\w+1},\aleph_\w)\twoheadrightarrow(\aleph_2,\aleph_1)$ fails.
\end{Theorem}

\subsection{Notation and definitions} 
The paper aims to be accessible and self-contained for any
reader with a strong background in combinatorial set theory. While no familiarity with forcing is required, 
a basic acquaintance with large cardinals combinatorics is assumed. 
When not otherwise explicitly stated \cite{jecST} is the standard source for notation
and definitions. For a regular cardinal $\theta$, we use $H(\theta)$
to denote the structure $\la H(\theta),\in,<\ra$ whose domain is the
collection of sets whose transitive closure is of size less than
$\theta$ and where $<$ is a predicate for a fixed well ordering of
$H(\theta)$. For cardinals $\lambda\leq\kappa$ we let
$[\kappa]^\lambda$ be the family of subsets of $\kappa$ of size
$\lambda$. In a similar fashion we define $[\kappa]^{<\lambda}$,
$[\kappa]^{\leq\lambda}$, $[X]^\lambda$, where $X$ is an arbitrary
set. If $X$ is an uncountable set and $\theta$ a regular cardinal,
$\E\subseteq[X]^\theta$ is unbounded if for every $Z\in [X]^\theta$,
there is $Y\in\E$ containing $Z$. $\E$ is bounded otherwise. For a
set of ordinals $X$, $\overline{X}$ denotes the topological closure
of $X$ in the order topology. For regular cardinals
$\lambda<\nu$, $S^\nu_{\lambda}$ denotes the subset of
$\nu$ of points of cofinality $\lambda$. In a similar fashion
we define $S^\nu_{<\lambda}$, $S^\nu_{>\lambda}$, etc... 
For the ease of the reader we will let $\theta<\lambda<\nu$ 
range over regular cardinals and $\kappa$ range over singular cardinals in most cases of cofinality $\theta$, moreover
unless otherwise stated the reader may safely assume that $\nu=\kappa^+$.
We say that a family $\mathcal{D}$ is covered by a family $\mathcal{E}$ if for
every $X\in\mathcal{D}$ there is a $Y\in\mathcal{E}$ such that
$X\subseteq Y$.

%%%%%%%%%%%%%%%%%%%%%%%%%%%%%%%%%%%%%%%%%%%%%%%%%%%%%%%%%%%%%%%%%%%%%%%%%%%%%%%%%%%%%%%%%%%%%%%%%%%%%%

\section{Covering matrices and the approachability ideal}
 
Shelah provides a characterization of the ideal $\I[\kappa^+,\kappa]$ which is suitable for our analysis. 
Let $\kappa$ be singular and let:
$$
d:[\kappa^+]^2\rightarrow\cf(\kappa).
$$
\begin{itemize}
\item $d$ is normal if $D(i,\beta)=\{\a<\beta:d(\a,\beta)\leq i\}$ has size less
than $\kappa$ for all $i$ and $\beta$, 

%\item $d$ is transitive if
%$d(\a,\gamma)\leq\max(d(\a,\beta),d(\beta,\gamma))$ for all $\a\leq\beta\leq\gamma$,

\item $\delta$ is $d$-approachable if there is $H$ unbounded in $\delta$ such
that $d[[H]^2]$ is bounded in $\cf(\kappa)$.

\end{itemize}

\noindent The following is an equivalent definition of
$\I[\kappa^+,\kappa]$ (theorem 3.28 \cite{eisHST}):

\begin{property} \label{shAIpr} Let $\kappa$ be singular of cofinality $\theta$.
$S\in\I[\kappa^+,\kappa]$ if and only if there are a normal coloring $d$ and a club $C\subseteq\kappa^+$ such that
$\delta$ is $d$-approachable for all $\delta\in S\cap C$.
\end{property}

\begin{proof}
We prove only the backward direction which is the one that we need. So assume $X$ is a subset of $\kappa^+$ such that for some normal
$d$ any $\delta\in X$ is $d$-approachable. Let $D(i,\beta)=\{\a<\beta:d(\a,\beta)\leq i\}$. We want to define a family
$\E=\{e_\a:\a<\kappa^+\}$ such that every point in $X$ is weakly approachable with respect to $\E$. To this aim fix a bijection $\phi:\theta\rightarrow\theta^2$ and let $\pi_0$ and $\pi_1$ be the projection maps of $\theta^2$ onto $\theta$. Notice that every ordinal $\delta$ below $\kappa^+$ can be decomposed uniquely as the sum $\delta=\alpha+i$ where $i<\theta$ and $\alpha$ is divisible by $\theta$. Now for every $\alpha<\kappa^+$ divisible by $\theta$ and for every $i<\theta$ set $e_{\alpha+i}=D(\pi_0\circ\phi(i),\alpha+\pi_1\circ\phi(i))$.
It is not hard to check that if $\delta$ is $d$-approachable, then it is weakly approachable with respect to $\E$.
\end{proof}

The coloring $d$ is determined by the matrix 
$\D(d)=\{D(i,\beta):i<\cf(\kappa),\beta<\kappa^+\}$ where $D(i,\beta)=\{\a<\beta:d(\a,\beta)\leq i\}$. It will be convenient for us to treat such matrices instead that the related coloring. Our aim is to show that mild reflection properties of a regular $\lambda<\kappa$ entail that for a suitably chosen normal coloring $d$
all points in $\kappa^+$ of cofinality $\lambda$ are $d$-approachable. This leads us to introduce and analyze the notion of a covering matrix.

\subsection{Covering matrices}
The reader is referred to \cite{viaMRL} for a detailed account of the results that are mentioned here without proof.

\begin{definition} \label{defCP}
For regular cardinals $\theta<\lambda$, $\D=\{D(i,\beta):i<\theta,\, \beta<
\lambda\}$ is a
$\theta$-covering matrix for $\lambda$ if:

\begin{itemize}
\item[\it (i)] $\beta=\bigcup_{i<\theta} D(i,\beta)$ for all $\beta$,
\item[\it (ii)] $D(i,\beta)\subseteq D(j,\beta)$ for all $\beta<\lambda$ and for
all $i<j<\theta$,
\item[\it (iii)]  for all $\beta<\gamma<\lambda$ and for all $i<\theta$, there is
$j<\theta$ such that $D(i,\beta)\subseteq D(j,\gamma)$.
\end{itemize}

\noindent A $\theta$-covering matrix $\D$ is transitive if $\a\in D(i,\beta)$
implies $D(i,\a)\subseteq D(i,\beta)$.

\smallskip

\noindent A $\theta$-covering matrix $\D$ is closed if $\sup X\in D(i,\beta)$ for all $X\in [D(i,\beta)]^{\leq\theta}$.

\smallskip

\noindent A $\theta$-covering matrix $\D$ is uniform if for all $\beta<\lambda$, $D(i,\beta)$ contains a club subset of $\beta$ for eventually all $i<\theta$.

\smallskip

\noindent $\beta_\D\leq\lambda$ is the least $\beta$ such that for
all $i$ and $\gamma$, $\otp(D(i,\gamma))<\beta$.
$\D$ is normal if $\beta_\D<\lambda$.

\end{definition}

\begin{example}
$d:[\kappa^+]^2\rightarrow\cf(\kappa)$ is normal if $\D(d)$ is a normal $\cf(\kappa)$-covering matrix on $\kappa^+$ with $\beta_\D=\kappa$.
\end{example}

We will be interested in the matrices produced by the following lemma:

\begin{lemma} \label{lemCM} 
For every singular cardinal $\kappa$, there is a uniform, closed, transitive $\cf(\kappa)$-covering matrix 
$\D$ on $\kappa^+$ with $\beta_\D=\kappa$.
\end{lemma}
\begin{proof} Let $\kappa$ be singular of cofinality $\theta$.
Fix $\{\kappa_i: i<\theta\}$ increasing sequence of regular cardinals converging to $\kappa$. Let 
$\phi_\alpha:\kappa\rightarrow\alpha$ be a surjection for all $\a<\kappa^+$ such that $\phi_\a[\kappa_i]$ contains a club subset of $\a$ whenever $\a$ is limit of cofinality smaller than $\kappa_i$. Now set
$D_0(i,\beta)=\phi_\beta[\kappa_i]$ for all $i<\theta$ and $\beta<\kappa^+$. Define by recursion over $\xi\leq\theta^+$ and limit and $n<\omega$:
\begin{itemize}
\item $D_{\xi+2n+1}(i,\beta)=\overline{D_{\xi+2n}(i,\beta)}$,
\item $D_{\xi+2n+2}(i,\beta)=\bigcup\{D_{\xi+2n+1}(i,\a):\a\in D_{\xi+2n+1}(i,\beta)\}$,
\item $D_{\xi}(i,\beta)=\bigcup\{D_{\eta}(i,\beta):\eta<\xi\}$.
\end{itemize}
Now set $D(i,\beta)=D_{\theta^+}(i,\beta)$ and check that $\D=\{D(i,\beta):i<\theta,\beta<\kappa^+\}$ is
a uniform, closed, transitive $\cf(\kappa)$-covering matrix 
$\D$ on $\kappa^+$ with $\beta_\D=\kappa$.
\end{proof}

\begin{definition} Let $\D=\{D(j,\beta):j<\theta,\beta<\lambda\}$ 
be a $\theta$-covering matrix on $\lambda$.

\smallskip

$\CP(\D)$ holds if there is $A$ unbounded subset
of $\lambda$ such that $[A]^\theta$ is covered by $\D$.

\smallskip

$S(\D)$ holds if there is $S$ stationary subet of $\lambda$ such that for all
families $\{S_i:i<\theta\}$ of stationary subsets of $S$ there are $j<\theta$ and $\beta<\lambda$ such that $S_i\cap
D(j,\beta)$ is non-empty for all $i<\theta$.
\end{definition}

We will come back to the relation between approachability and covering matrices at the end of this section, we now aim
to investigate the consistency of $S(\D)$ and $\CP(\D)$ for a large variety of covering matrices $\D$.

\subsection{Consistency of $\CP(\D)$ and $S(\D)$}

\begin{fact} \label{faconSD}Assume $R(\lambda,\theta)$ holds and $\D$ is a uniform $\theta$-covering matrix on $\lambda$.
Then $S(\D)$ holds. 
\end{fact}

\begin{proof} Let $\D$ be a uniform $\theta$-covering matrix on $\lambda$
and $\{S_i:i<\theta\}$
be a family of stationary subsets of $S$. By $R(\lambda,\theta)$ find $\delta$ such
that $S_i$ reflects on $\delta$ for all $i<\theta$. Now $\D$ is uniform, so there is a $j<\theta$ such that $D(j,\delta)$ contains a club subset of $\delta$. Thus $S_i\cap D(j,\delta)$ is non-empty for all $i<\theta$. Since the family
$\{S_i:i<\theta\}$ is arbitrary $S(\D)$ holds as witnessed by $S$.
\end{proof}

\begin{corollary} $\MM$ implies $S(\D)$ for all uniform $\theta$-covering matrices $\D$ on $\lambda$ whenever
$\lambda>\aleph_1$ is a regular cardinal and $\aleph_1\geq\theta$.
\end{corollary}
\begin{proof} $\MM$ implies $R(\lambda,\aleph_1)$ holds as witnessed by $S^\lambda_\w$ for all regular $\lambda>\aleph_1$.
\end{proof}

In \cite{viaMRL} it is shown the following:

\begin{theorem}\label{thPFACP} 
$\PFA$ implies $\CP(\D)$ for all $\w$-covering matrix $\D$ on a regular
$\lambda>\aleph_2$.
\end{theorem}

We now investigate the relation between $\CP(\D)$ and $S(\D)$ and show that they are
equivalent whenever $\D$ is transitive and closed.

\subsection{When are $\CP(\D)$ and $S(\D)$ equivalent?}

\begin{proposition}\label{lemCPS}
Let $\D$ be a $\theta$-covering matrix on $\lambda$.
The following holds:
\begin{description}
\item[\it (i)] $\CP(\D)$ implies $S(\D)$ whenever $\D$ is closed,
\item[\it (ii)] $S(\D)$ implies $\CP(\D)$ whenever $\D$ is transitive.
\end{description}
\end{proposition}
\begin{proof} We first show {\it (i)}. We will actually show that
$\CP(\D)$ implies
$S(\D)$ is witnessed by $S^\lambda_\theta$.
So let $\{S_i:i<\theta\}$ be a family of stationary subsets of $S^\lambda_\theta$. By
$\CP(\D)$, there is $X$
unbounded in $\lambda$ such that $[X]^\theta$ is covered by $\D$.
We claim that
$[\overline{X}\cap S^{\lambda}_{\theta}]^\theta$ is covered
by $\D$. To see this, let $Z$ be in this latter set and find
$Y\subseteq X$ of size $\theta$ such that $Z\subseteq\overline{Y}$.
Now find $i$ and $\beta$ such that $Y\subseteq D(i,\beta)$. Since
$D(i,\beta)$ is closed under sequences of size at most $\theta$, $Z\subseteq\overline{Y}\subseteq
D(i,\beta)$.

Now pick $M\prec H(\lambda)$ with $\lambda$ large enough
such that $|M|=\theta\subseteq M$ and $\theta, X,
\{S_i:i<\theta\}\in M$. Now $S_i\cap\overline{X}$ is non-empty for all $i<\theta$. By
elementarity, $M$ sees this and so $M\cap S_i\cap\overline{X}$ is non-empty for all
$i<\theta$. However
$M\cap\overline{X}\cap S^{\lambda}_{\theta}$ has size $\theta$
so there are $j$ and $\beta$ such that $M\cap\overline{X}\cap
S_{\theta}^{\lambda}\subseteq D(j,\beta)$. So $S_i\cap
D(j,\beta)$ is non-empty for all $i<\theta$. This proves the first implication.

\smallskip 

We now show {\it (ii)}. So assume $S(\D)$ holds for
a transitive $\theta$-covering matrix $\D$ on $\lambda$. Let
$S$ witness $S(\lambda,\theta)$ and $T_i$ be the set of $\a\in S$ such that
$$
S_\a^i=\{\beta\in S\setminus\a:\a\in D(i,\beta)\}
$$ 
is stationary. It is
straightforward to see that for some $i<\theta$, $T_i$ is stationary. We aim to show that $[T_i]^\theta$ is covered by $\D$: let $X\in[T_i]^\theta$ and consider the family of stationary sets
$\{S_\a^i:\a\in X\}$. Since $X$ has size $\theta$, by $S(\D)$ there are some
$j<\theta$ and $\delta<\lambda$ such that $S_\a^i\cap D(j,\delta)$ is non-empty for all $\a\in X$. W.l.og. we can suppose that $j\geq i$.
Now for any $\a\in X\subseteq T_i$, there is $\beta_\a\in
D(j,\delta)\cap S_\a^i$, i.e. $\beta_\a$ is such that $\a\in D(i,\beta_\a)$. Since
$\D$ is a transitive covering matrix and $j\geq i$, 
$$
\a\in D(i,\beta_\a)\subseteq D(j,\beta_\a)\subseteq D(j,\delta).
$$ 
This means that $X\subseteq D(j,\delta)$. Since $X$ is arbitrary we can conclude that
$[T_i]^{\theta}$ is covered by $\D$.
\end{proof}

\subsection{A weak form of diagonal reflection}\label{subsecWDR}
We aim to show that $\CP(\D)$ or $S(\D)$ strongly limits the kind of behavior a $\theta$-covering matrix $\D$ on $\lambda$ may have. We shall now see that $\CP(\D)$ plus suitable assumptions on the proportion between the width $\theta$ and the height $\lambda$ of $\D$ imply that there is an unbounded subset of $\lambda$ such that all its initial segments are covered by $\D$. Once this is achieved, it will be easy to conclude that $R(\lambda,\theta)$ implies that all points of cofinality $\lambda$ below $\kappa^+$ are weakly approachable whenever $\kappa>\lambda$ is a singular cardinal of cofinality $\theta$. We now prove that a weak form of diagonal reflection of stationary sets on many covering matrices $\D$ follows from $S(\D)$ or $\CP(\D)$.

\begin{lemma}\label{lemappr} 
Assume $\D$ is a $\theta$-covering matrix on $\lambda$, $S(\D)$ holds as witnessed by $S$ and that either $\lambda$ is $\theta$-inaccessible or $\lambda\in(\theta,\theta^{+\w})$. Then
for all families $\{S_\beta:\beta<\lambda\}$ of stationary subsets of $S$ there are $\delta<\lambda$
and $i<\theta$ such that $S_\a\cap D(i,\delta)$ is non-empty for all $\a<\delta$.
\end{lemma}

\begin{lemma}\label{lemapprbis} 
Assume $\D$ is a $\theta$-covering matrix on $\lambda$. $\CP(\D)$ holds as witnessed by $T$ and that either $\lambda$ is $\theta$-inaccessible or $\lambda\in(\theta,\theta^{+\w})$. Then there are stationarily many $\delta<\lambda$ such that $T\cap\delta\subseteq D(i,\delta)$ for some $i<\theta$.
\end{lemma}

We give a detailed proof of the first lemma. The second lemma is proved by a self evident step by step modification of this argument. 

\begin{fact} \label{farefl1}
Let $\theta<\lambda<\nu$ be regular cardinals such that
$\lambda^{\theta}<\nu$, $\D=\{D(j,\beta):j<\theta,\beta<\nu\}$ a 
$\theta$-covering matrix on $\nu$ and assume $S(\D)$ holds as witnessed by $S$. Let
$\{S_i:i<\lambda\}$ be a family of stationary subsets of $S$. Then there are
$j<\theta$ and $\beta<\nu$ such that $S_i\cap D(j,\beta)$ is non-empty for all
$i<\lambda$.
\end{fact}

\begin{proof} Assume not and let $\{S_i:i<\lambda\}$ contradict the fact.
We need to find
$j<\theta$ and $\beta<\nu$ such that $S_i\cap D(j,\beta)$ is non empty for all
$i<\lambda$.
For $X\in [\lambda]^\theta$ let by $S(\D)$, $k_X<\theta$ and
$\beta_X<\nu$ be such that $S_i\cap
D(k_X,\beta_X)$ is non-empty for all $i\in X$. By our assumptions,
$\lambda^{\theta}<\nu$. For
this reason $\beta=\sup_{X\in [\lambda]^\theta}\beta_X<\nu$. Now by
property {\it (ii)} of $\mathcal{D}$, we have that for all $X\in
[\lambda]^\theta$, $D(k_X,\beta_X)\subseteq D(j_X,\beta)$ for some $j_X<\theta$.
Let $\C_j$
be the set of $X$ such that $j_X=j$. Now notice that for at least
one $j$, $\C_j$ must be unbounded in $[\lambda]^\theta$, otherwise
$[\lambda]^\theta$ would be the union of $\theta$-many bounded subsets, which is
not possible since $\lambda$ has cofinality different from $\theta$. Then
$S_i\cap D(j,\beta)$ is non-empty for all $i<\lambda$, since every $i\in \lambda$ is
in some
$X\in\C_j$, as $\C_j$ is unbounded. This completes the proof
of the fact.
\end{proof}

\begin{fact} \label{farefl2}
Assume $\lambda\in(\theta,\theta^{+\w})$, $\nu>\lambda$ is regular and $S(\D)$ holds for some $\theta$-covering matrix $\D$ on $\nu$ and is witnessed by
$S$. Let $\{S_i:i<\lambda\}$ be a family of stationary subsets of $S$. Then there are
$j<\theta$ and $\beta<\nu$ such that $S_i\cap D(j,\beta)$ is non-empty for all
$i<\lambda$.
\end{fact}

\begin{proof} Proceed by induction on $n$ so assume the claim holds for
$\theta^{+n}$ and let $\lambda=\theta^{+(n+1)}$ and $\{S_i:i<\lambda\}$ be a family of
stationary subsets of $S$.
By the inductive assumption for all $i<\lambda$, there are $k_i<\theta$ and
$\beta_i<\nu$ such that $S_j\cap D(k_i,\beta_i)$ is non-empty for all $j<i$.
Since $\lambda<\nu$ there is $\beta<\nu$ larger than all $\beta_i$. Now by
property {\it (ii)} of $\D$ we have that for all $i<\lambda$ there is $j_i<\theta$
such that $D(k_i,\beta_i)\subseteq D(j_i,\beta)$. Find $U$ unbounded subset of
$\lambda$ such that for all $i\in U$, $j_i=j$. We can conclude that $S_l\cap
D(j,\beta)$ is non-empty for all $l<\lambda$, since
$S_l\cap D(k_i,\beta_i)$ is non-empty provided $l<i$ and $i\in U$ and
$D(k_i,\beta_i)\subseteq D(j,\beta)$.
\end{proof}

We are now ready to prove lemma \ref{lemappr}.

\begin{proof} Assume not and let $\{S_\beta:\beta<\nu\}$ contradict the
lemma. 
For each $\delta$ of cofinality larger then $\theta$, let $\gamma_\delta<\delta$ be
the least such that for all $i<\theta$ there is $\gamma_i^\delta<\gamma_\delta$ such
that $S_{\gamma_i^\delta}\cap D(i,\delta)$ is empty. Find $A$ stationary subset of
$\nu$ such that $\gamma_\delta=\gamma$ for all $\delta\in A$. By our assumption
on $\nu$ and facts \ref{farefl1} and \ref{farefl2}, we know that there are
$i<\theta$ and $\delta_0<\nu$ such that $S_\a\cap D(i,\delta_0)$ is non empty for
all $\a<\gamma$. Pick $\delta\in A\setminus\delta_0$ and $j<\theta$ such that
$D(i,\delta_0)\subseteq D(j,\delta)$. Then we get that $S_{\gamma_j^\delta}\cap
D(j,\delta)$ is non-empty since $S_{\gamma_j^\delta}\cap D(i,\delta_0)$ is non-empty
and $D(i,\delta_0)\subseteq D(j,\delta)$. This contradicts the very definition of
$\gamma_j^\delta$.
\end{proof}

\noindent In particular we have shown the following:

\begin{fact} \label{keyfa}Assume $\lambda$ is either $\theta$-inaccessible or $\lambda\in(\theta,\theta^{+\w})$ and $S(\D)$ holds for a transitive $\theta$-covering matrix $\D$ on $\lambda$. Then there is $A$ unbounded subset of $\lambda$ such that $[A]^{<\lambda}$ is covered by $\D$.
\end{fact}

\subsection{Main result} We are now in the position to state our main result:
\begin{theorem} \label{mainth2}
Assume $\kappa$ is singular of cofinality $\theta$ and a regular $\lambda<\kappa$ is either 
$\theta$-inaccessible or in $(\theta,\theta^{+\w})$ and such that $S(\D)$ (or equivalently $\CP(\D)$) holds for all uniform, closed and transitive $\theta$-covering matrices $\D$ on $\lambda$. Then club many points in $\kappa^+$ of cofinality $\lambda$ are approachable.
\end{theorem}

\begin{proof} Fix $d:[\kappa^+]^2\rightarrow\theta$ such that $D(d)=\{D(i,\beta):i<\theta,\beta<\kappa^+\}$ is a normal, uniform, closed and transitive
$\theta$-covering matrix on $\kappa^+$, where $D(i,\beta)=\{\alpha<\beta:d(\alpha,\beta)\leq i\}$. Such a $d$ exists by lemma \ref{lemCM}. By property
\ref{shAIpr} it is enough to show that all points of cofinality $\lambda$ are $d$-approachable. Let $\beta$ be such that $\cf(\beta)=\lambda$.
Find $A=\{\delta_\xi:\xi<\lambda\}$ closed and unbounded subset of $\beta$ of minimal order-type. Let $\pi$
be the transitive collapse of $A$ on $\lambda$ and $\E=\{E(i,\xi):i<\theta,
\xi<\lambda\}$ be the matrix whose entries are the sets $\pi[D(i,\delta_\xi)\cap
A]$. Then $\E$ is a uniform, closed and transitive $\theta$-covering matrix on $\lambda$. 
By $S(\E)$ and fact \ref{keyfa}, there is $B$ unbounded subset of $\lambda$ such that $[B]^{<\lambda}$ is covered by
$\E$. Thus $B\cap\eta\subseteq E(i_\eta,\xi_\eta)$ for some
$i_\eta<\theta$ and $\xi_\eta\in B\setminus\eta$ for all $\eta\in B$. Refine $\{\xi_\eta:\eta\in B\}$ to un unbounded subset $C$ such that $\xi_\eta<\gamma$ for all $\xi_\eta<\xi_\gamma\in C$. Thus $\xi_\eta\in B\cap\gamma\subseteq E(i_\gamma,\xi_\gamma)$ for all $\xi_\eta<\xi_\gamma\in C$. Let $D$ be an unbounded subset of $C$ such that for some fixed $j$, $i_\eta=j$ for all $\xi_\eta\in D$. Now if $\xi_\eta<\xi_\gamma\in D$ we have that $\xi_\eta\in B\cap\gamma\subseteq E(j,\xi_\gamma)$ i.e.
$d(\pi^-1(\xi_\eta),\pi^{-1}(\xi_\gamma))\leq j$, i.e $\pi^{-1}[D]$ witnesses that
$\beta$ is $d$-approachable.
\end{proof}

\section{Joint reflection of stationary sets}\label{secRP}

We briefly analyze the consistency strength of the hypothesis of theorem \ref{mainth}.

\begin{fact} \label{fareflWC}
$R^*(\lambda)$ holds if $\lambda$ is weakly compact.
\end{fact}

\begin{proof} Recall the following characterization of weak compactness: $\lambda$ is
weakly compact
if for every transitive model $M$ of $\ZFC$ minus the powerset axiom such that $M$
has size $\lambda$ and $H(\lambda)\subseteq M$, there is an elementary embedding of
$M$ into a transitive structure $N$ with critical point $\lambda$. Now let
$\{S_i:i<\lambda\}$ be any family of stationary subsets of $\lambda$.
To prove $R^*(\lambda)$ we must find a $\delta<\lambda$ such that $S_\a$ reflects on
$\delta$ for all $\a<\delta$.
Let $M$ be a structure as above such that $\{S_i:i<\lambda\}\in M$.
Let $j:M\rightarrow N$ be elementary with $N$ transitive and critical point of
$j=\lambda$.
Then $j(\{S_\a:\a<\lambda\})=\{T_\a:\a<j(\lambda)\}$ and $j(S_\a)\cap\lambda=S_\a$ for
all $\a<\lambda$.
Thus $N$ models that there is $\delta<j(\lambda)$ (namely $\delta=\lambda$) such that
for all $\a<\delta$,
$T_\a$ reflects on $\delta$. By elementarity of $j$ there is $\delta<\lambda$ such
that $S_\a$ reflects on $\delta$ for all $\a<\delta$ and we are done.
\end{proof}

\noindent Larson (unpublished) has proved that $\MM$ implies $R^*(\aleph_2)$ while
it is apparent already in the paper of Foreman, Magidor and Shelah \cite{formagsheMM} 
that $\MM$ implies $R(\lambda,\aleph_1)$ for all regular $\lambda>\aleph_1$. On the other hand Magidor
\cite{magSW2} has shown that $R^*(\aleph_n)$ is equiconsistent with $\aleph_n$ being
weakly compact in $L$. Notice however that a model of $R(\aleph_n,\aleph_0)$ and
$R(\aleph_{n+1},\aleph_0)$ subsume already very large cardinal assumptions since it can be seen that
$R(\aleph_n,\aleph_0)$ implies failure of\footnote{For example Jensen produces from $\Box(\lambda)$ a transitive and closed $\omega$-covering matrix $\D$ on $\lambda$ such that $\CP(\D)$ cannot hold. A proof of this result by Todor\v{c}evi\'c can be found in Todor\v{c}evi\'c's book \cite{todWOO} or in section 2.2.1 of \cite{viaPHD}.} $\Box(\aleph_n)$ and Schimmerling has shown
that failure of $\Box(\aleph_n)$ for two consecutive cardinals implies projective deteminacy
\cite{schPD}.
Another scenario suggested by Foreman to obtain $R^*(\lambda)$ is the following:
\begin{lemma}
Assume that $\I$ is a $\lambda$-complete,
fine ideal which concentrates on $[\kappa]^{<\lambda}$ for some $\kappa\geq\lambda$ 
and such that $P_\I=P([\kappa]^{<\lambda})/\I$ is a proper forcing. Then $R^*(\lambda)$ holds.
\end{lemma}

\begin{proof} First of all $\I$ is precipitous since $P_\I$
is proper (\cite{forHST} Proposition 4.10). Let $G$ be a generic filter for $P_\I$. Then the
ultrapower $M=V^{([\kappa]^{<\lambda})}\cap V/G$ defined in $V[G]$ is well-founded.
Let $j:V\rightarrow M$ be the associated generic elementary embedding. Since $\I$ is
$\lambda$-complete and fine, we have that the critical point of $j$ is $\lambda$.
Now let $\{S_\a:\a<\lambda\}\in V$ be a family of stationary subsets of
$S^\lambda_{\aleph_0}$. It is clear that
$M$ models that $j(\{S_\a:\a<\lambda\})=\{T_\a:\a<j(\lambda)\}$ is a family of
stationary subset of $S^{j(\lambda)}_{\aleph_0}$. Now
$T_\a=j(S_\a)$ and $j(S_\a)\cap\lambda=S_\a$ for all $\a<\lambda$. Since
$P$ is proper, $S_\a$ is a stationary subset of $\lambda$ in $V[G]$ so it is
certainly a stationary subset of $\lambda$ in $M$.
Then $M$ models that $j(S_\a)$ reflects on $\lambda$ for all $\a<\lambda$. Now the
argument to show that $S^*(\lambda)$ holds in $V$ is as in fact \ref{fareflWC}.
\end{proof}

Notice that we've hidden a large cardinal assumption
in the requirement that $P$ is proper. The hypothesis of the lemma are satisfied by the non-stationary ideal on $\aleph_2$ in the generic extension by a Levy collapse of a measurable $\lambda$ to $\aleph_2$. In this case the quotient algebra is even countably complete. \cite{forHST} is a survey on
generic large cardinals. We now turn to an application of the main theorems \ref{mainth} and \ref{mainth2} to Chang conjectures.

%%%%%%%%%%%%%%%%%%%%%%%%%%%%%%%%%%%%%%%%%%%%%%%%%%%%%%%%%%%%%%%%%%%%%%%%%%%%%%%%%%%%%%%%%%%%%%%%%%%%%%%%%%%%%%%%%%

\section{$R(\aleph_2,\aleph_0)$ denies $(\aleph_{\omega+1},\aleph_\omega)\twoheadrightarrow(\aleph_2,\aleph_1)$} \label{secCC}

Recall that the Chang conjecture $(\lambda,\kappa)\twoheadrightarrow(\theta,\nu)$ holds for $\lambda>\kappa\geq\theta>\nu$ if for every structure $\langle Y,\lambda,\kappa,...\rangle$ with predicates for $\lambda$ and $\kappa$ there is $X\prec Y$ such that $|X\cap\lambda|=\theta$ and $|X\cap\kappa|=\nu$. We will also be interested in the principles of the form $(\kappa,\lambda)\twoheadrightarrow(\theta,<\theta)$ which are likewise defined. It is known that $(\aleph_2,\aleph_1)\twoheadrightarrow(\aleph_1,\aleph_0)$ as well as many other Chang conjectures are consistent relative to appropriate large cardinals assumptions. For example it is possible to see
that $(j(\kappa^{+\theta}),j(\kappa^{+\gamma}))\twoheadrightarrow(\kappa^{+\theta},\kappa^{+\gamma})$ whenever $\kappa$ is the critical point of a 2-huge embedding and $\gamma<\theta<\kappa$. Developing on this, Levinsky, Magidor and Shelah in \cite{levmagsheCCW} showed that $(\aleph_{\w+1},\aleph_\w)\twoheadrightarrow(\aleph_1,\aleph_0)$ is consistent relative to the existence of a 2-huge cardinal. However all the known examples of a consistent $(\kappa^+,\kappa)\twoheadrightarrow(\theta^+,\theta)$ where $\kappa$ is singular and $\theta$ regular are such that $\theta=\cf(\kappa)$. Thus a folklore problem in this field is the following:

\begin{problem} \label{mainpb}
Is it consistent that $(\kappa^+,\kappa)\twoheadrightarrow(\theta^+,\theta)$ for some regular $\theta$ and singular $\kappa$ of cofinality smaller than $\theta$?
\end{problem}

First of all it is a simple fact that such Chang conjectures affect cardinal arithmetic:

\begin{fact} \label{faca}Assume $(\kappa^+,\kappa)\rightarrow(\theta^+,\theta)$ for some singular $\kappa$. Then $\theta^+\leq\theta^{\cf(\kappa)}$.
\end{fact}

\begin{proof} Notice that $\kappa^{\cf(\kappa)}>\kappa$. Now assume $(\kappa^+,\kappa)\rightarrow(\theta^+,\theta)$. Fix $\lambda>\kappa^+$ regular and large enough and let $H(\lambda)$ denotes the family of sets whose transitive closure has size less than $\theta$. Fix $M\prec \langle H(\lambda),\kappa^+,\kappa,....\rangle$ with $|M\cap\kappa^+|=\theta^+$ and $|M\cap\kappa|=\theta$. Pick a family $\{X_\a:\a<\kappa^+\}\in M$ of distinct elements of $[\kappa]^{\cf(\kappa)}$. By elementarity of $M$, $X_\a\cap M\neq X_\beta\cap M$ for all $\a,\beta\in M\cap\kappa^+$. Thus we have a family of $\theta^+$ distinct elements of $[M\cap\kappa]^{M\cap\cf(\kappa)}$. Now $|M\cap\kappa|=\theta$ and $|M\cap\cf(\kappa)|\leq\cf(\kappa)$. Thus $\theta^+\leq|[M\cap\kappa]^{M\cap\cf(\kappa)}|\leq\theta^{\cf(\kappa)}$.
\end{proof}

Cummings in \cite{cumCCW} has shown that these Chang conjectures can be studied by means of $\pcf$-theory and has obtained several other restrictions on the combinatorics of the singular cardinals $\kappa$ which may satisfy an instance of the above problem.
For example he has shown that these Chang conjectures subsume the existence of very strong large cardinals, i.e. out of the scope of analysis of the current inner model theory: it can be argued by the analysis brought up in \cite{cumCCW} that  $(\aleph_{\w+1},\aleph_\w)\twoheadrightarrow(\aleph_n,\aleph_{n-1})$, then $\Box_{\aleph_\w}$ fails and $\SCH$ holds at $\aleph_\w$. Moreover a result by Shelah shows that $n$ cannot be greater than\footnote{Cummings' analysis relies on the notion of good (or in Kojamn's notation \cite{kojEUB} flat) points for a scale on $\prod_n\aleph_n$ and his main observation (Lemma 3.1 of \cite{cumCCW}) is that if $(\aleph_{\w+1},\aleph_\w)\twoheadrightarrow(\aleph_n,\aleph_{n-1})$ holds, then there are stationarily many non-good points of cofinality $\aleph_n$. On the other hand Shelah has shown that club many points of cofinality $\aleph_n$ are good (or flat) for a scale on $\prod_n\aleph_n$ if either $\aleph_n>2^{\aleph_0}$ (Exercise 2.9-2, Lemma 2.12 and Theorem 2.13 of \cite{abrmagPCFHST}) or $n>3$ (Theorem 2.13 and Lemmas 2.12 and 2.19 of \cite{abrmagPCFHST}).} $3$.   
We can decrease $3$ down to $1$ and greatly simplify their argument avoiding any mention of scales in the case that $R(\aleph_n,\aleph_0)$ holds:

\begin{theorem} Assume $R(\aleph_n,\aleph_0)$ holds for some $n>1$. 
Then $(\kappa^+,\kappa)\twoheadrightarrow(\aleph_n,\aleph_{n-1})$ fails for all singular $\kappa$ of countable cofinality.
\end{theorem}
\begin{proof} Towards a contradiction let $\kappa$ and $n$ be counterexamples to the theorem. Fix $\mu>2^{\kappa^+}$ regular and large enough and $M\prec H(\mu)$ structure containing all relevant information and such that $|M\cap\kappa|<\aleph_n$ and $|M\cap\kappa^+|=\aleph_n$. First of all:
\begin{claim}
$\otp(M\cap\kappa^+)=\aleph_n$.
\end{claim} 
Suppose otherwise and let $\gamma\in M$ be such that $\otp(M\cap\gamma)=\aleph_n$. Then $\gamma\in M\cap(\kappa,\kappa^+)$. We claim that $M$ models that $\gamma$ is a regular cardinal, which gives the desired contradiction since by elementarity $\gamma$ would be in the universe a regular cardinal in $(\kappa,\kappa^+)$ which is impossible. So suppose $M$ models $\gamma$ is not a cardinal, then in $M$ there is a bijection $\phi$ of $\gamma$ onto $\kappa$. If we take the transitive collapse $\pi_M$ of $M$, $\pi_M(\phi)$ is now a bijection of $\pi_M(\gamma)=\aleph_n$ onto
$\pi_M(\kappa)$ which is an ordinal of size less than $\aleph_n$. Contradiction. \hfill$\square$\medskip

Fix in $M$ a transitive, uniform and closed $\omega$-covering matrix $\D$ on $\kappa^+$ with $\beta_\D=\kappa$.
Let $X=M\cap\kappa^+$ and $\delta_M=\otp(M\cap\kappa)$. 

We will use the following:
\begin{fact} $\otp(\overline{Y})\leq\otp(Y)+1$ for any set of ordinals $Y$. 
\end{fact}
This can be proved by induction on the ordertype of $Y$. \hfill$\square$\medskip

Consider now the sets $D(i,\beta)\cap\overline{X}$ for $i<\w$ and $\beta\in\overline{X}$.

\begin{claim} $D(i,\beta)\cap\overline{X}<\delta_M$ for all $i<\w$ and $\beta\in\overline{X}\cap\kappa^+$. 
\end{claim}
Notice that for all $\beta\in\overline{X}\cap\kappa^+$ and $i<\omega$, 
$D(i,\beta)\cap\overline{X}\subseteq D(j,\gamma)\cap\overline{X}$ for some $\gamma\in X\setminus\beta$ and $j\geq i$ since $\D$ is an $\w$-covering matrix. So it is enough to prove the claim for all $i<\w$ and $\beta\in X$. Now if $\beta\in X$, $D(i,\beta)\in M$. Since $\otp(D(i,\beta))<\kappa$, by elementarity of $M$, we get that $\otp(D(i,\beta)\cap X)<\otp(X\cap\kappa)=\delta_M$. Now $\delta_M$ is a limit ordinal, so $\otp(D(i,\beta)\cap\overline{X})\leq\otp(D(i,\beta)\cap X)+1<\delta_M$. The claim is now proved. \hfill$\square$\medskip

Let $\pi_{\overline{X}}$ be the transitive collapse of $\overline{X}$ and set 
$$
\E=\{E(i,\beta):i<\w,\beta<\aleph_n\}
$$ 
where $E(i,\beta)=\pi_{\overline{X}}[D(i,\pi_{\overline{X}}^{-1}(\beta))\cap\overline{X}]$.
Now $\E$ is a transitive, uniform, $\w$-covering matrix on $\aleph_n$ with $\beta_\E\leq\delta_M<\aleph_n$.
By $R(\aleph_n,\aleph_0)$, $S(\E)$ holds, so there is $A$ unbounded in $\aleph_n$ such that $[A]^{\aleph_0}$ is covered by $\E$.
By lemma \ref{lemapprbis} $[A]^{\aleph_{n-1}}$ is covered by $\E$. So find $\gamma$ such that $A\cap\gamma$ has order type bigger than $\delta_M$.
Now there are $j$ and $\xi$ such that $A\cap\gamma\subseteq E(j,\xi)$, then:
$$  
\delta_M<\otp(A\cap\gamma)\leq\otp(E(j,\xi))<\delta_M.
$$
This is the desired contradiction which proves the theorem.
\end{proof}

\section{Some open questions and some comments}

The original question by Magidor and Foreman \cite{formagVWS} remains open:

\begin{problem} 
Is it consistent that $S^{\aleph_{\w+1}}_{\aleph_2}\not\in\I[\aleph_{\w+1}]$?
\end{problem}

Foreman and Cummings have indipendently shown that it is possible to force rightaway in $\ZFC$ by a cardinal preserving forcing a transitive, uniform $\w$-covering matrix on $\aleph_2$ such that $S(\D)$ fails. Veli\v{c}kovi\'c noticed that it is possible to obtain further counterexamples to $S(\D)$ using Todor\v{c}evi\'c's techniques of minimal walks over a $\Box(\aleph_2)$-sequence\footnote{Todor\v{c}evi\'c's book \cite{todWOO} gives a complete account of the method of minimal walks and of its applications.}. So the reflection hypothesis on $\aleph_2$ are needed to obtain
that $S(\D)$ holds for all uniform and transitive $\w$-covering matrix $\D$ on $\aleph_2$. On the other hand no strategy to force $S^{\aleph_{\w+1}}_{\aleph_2}\not\in\I[\aleph_{\w+1}]$ seems currently available.

A negative answer to the above problem would entail also a negative answer to question \ref{mainpb}, i.e.:
\begin{problem}
Is it consistent that $(\aleph_{\omega+1},\aleph_\omega)\twoheadrightarrow(\aleph_2,\aleph_1)$?
\end{problem}
It seems more fruitful to attack this problem directly by means of Shelah's analysis of the existence of exact upper bounds for families of 
ordinal functions in $Ord^\w$ (see \cite{cumCCW} and section 2 of \cite{abrmagPCFHST}). For example using these techniques we can already prove that $(\aleph_{\omega+1},\aleph_\omega)\twoheadrightarrow(\aleph_n,\aleph_{n-1})$ fails if $n>3$.

A comment on our main theorem \ref{mainth} is in order: the theorem entails that in a model of $\MM$ $S^{\kappa^+}_\lambda\in\I[\kappa^+,\kappa]$ for all $\kappa$ of countable cofinality and for all regular $\lambda<\kappa$ which are $\w$-inaccessible i.e. which are not the successor of a cardinal of countable cofinality. We expect this to be close to the best possible result for models of $\MM$.
For example consider the following scenario: $\kappa$ is a supercompact cardinal and 
$(\lambda^{+\w+1},\lambda^{+\w})\twoheadrightarrow(\nu^{+\w+1},\nu^{+\w})$ for some $\nu>\kappa$ holds as witnessed by structures $M$ such that $M^\kappa\subseteq M$. This occurs if there is a $2$-huge cardinal larger than $\kappa$. Now force $\MM$ collapsing $\kappa$ to $\aleph_2$. In the resulting generic extension $\MM$ holds and the chain condition of the forcing is small enough to preserve the truth of $(\lambda^{+\w+1},\lambda^{+\w})\twoheadrightarrow(\nu^{+\w+1},\nu^{+\w})$. This Chang conjecture already implies that
$S^{\lambda^{+\w+1}}_{\nu^{+\w+1}}\not\in\I[\lambda^{+\w+1}]$.

\bibliographystyle{plain}
\bibliography{biblioMV}
\end{document}